\documentclass[a4paper,12pt]{article}
\usepackage{times, url}
\usepackage{amsmath,amssymb,amsthm,amsfonts}
\usepackage{graphicx}
\usepackage{dcolumn}
\usepackage{bm}
\usepackage{color}

\begin{document}
\numberwithin{equation}{section}
\vspace{10mm}

\begin{center}
{\LARGE \bf  Fast converging series for zeta numbers in terms of polynomial representations of Bernoulli numbers}
\vspace{8mm}

{\large \bf J. Braun$^1$, D. Romberger$^2$, H. J. Bentz$^3$}
\vspace{3mm}

$^1$Department Chemie, Ludwig-Maximilians-Universit\"at \\ 
M\"unchen, 81377 M\"unchen, Germany \\
e-mail: \url{juergen.braun@cup.uni-muenchen.de}
\vspace{2mm}

$^2$ Fakult\"at IV, Hochschule Hannover \\ 
Ricklinger Stadtweg 120, 30459  Hannover, Germany \\
e-mail: \url{detlef.romberger@hs-hannover.de}
\vspace{2mm}

$^3$ Institut f\"ur Mathematik und Angewandte Informatik \\ 
Samelsonplatz 1, 31141 Hildesheim, Germany \\
e-mail: \url{bentz@cs.uni-hildesheim.de}

\end{center}
\vspace{10mm}

\noindent
{\bf Abstract:} 
In this work we introduce a new polynomial representation of the Bernoulli numbers in terms of polynomial sums
allowing on the one hand a more detailed understanding of their mathematical structure and on the other
hand provides a computation of $B_{2n}$ as a function of B$_{2n-2}$ only. Furthermore, we show that
a direct computation of the Riemann zeta-function and their derivatives at k $\in \mathbb Z$ is possible
in terms of these polynomial representation. As an explicit example, our polynomial Bernoulli number
representation is applied to fast approximate computations of $\zeta$(3), $\zeta$(5) and $\zeta$(7). 

\vspace{10mm}
\section{Introduction} 
Much was written on Bernoulli numbers \cite{Sal1898,Comt74,Brass85,Apostol98,Hu13} since Jacob Bernoulli
has discovered these fascinating quantities \cite{Arscon17}. Their importance is widely documented as
they appear in a variety of scientific fields of importance \cite{Kum1850,Brass85,Diu89,Ara01}. In particular,
their close relationship with the Riemann zeta-function \cite{Sita86,Bentz94,Adam98}, shows the need for
a better understanding of their mathematical structure. Bernoulli numbers are rational, their fractional
part is known by the results of Karl von Staudt and Thomas Clausen, often cited as the Staudt-Clausen theorem
\cite{Ada1872,Cal07}. Furthermore, several direct representations exist \cite{Comt74,Gould72,Dil89,Todo91},
all show up as complicated expressions in form of double-sums providing not much information on the enumerator
as well as on the denominator of these numbers. No simple rules to compute Bernoulli numbers had been established
so far and in consequence explicit formulas are more or less of pure academic interest \cite{Brass85}. On the
other hand simpler representations could be of great practical importance, for example, in the field of analytical
number theory. For a direct computation of the Riemann zeta-function and their derivatives for all k $\in \mathbb Z$
in terms of Bernoulli numbers is possible. Let us remark first on these relations. Some of such relations are
history know like the famous Euler representation \cite{Hecke44,Hav03}
\begin{eqnarray}
\zeta (2n) = (-)^{n+1} B_{2n} \frac{(2\pi)^{2n}}{2(2n)!}~,
\end{eqnarray}
or for negative integers \cite{Adam98}
\begin{eqnarray}
\zeta (-n) = -  \frac{B_{n+1}}{n+1}~.
\end{eqnarray}
More delicate is the direct computation of $\zeta (2n+1)$ and of $\zeta'(2n)$. Explicit relations
were found by us and published in 1994 \cite{Bentz94} and similar one later by Adamchik et al. \cite{Adam98,Lima14}.
Some remarks on the corresponding computational scheme will be given with an explicit example later on,
where we demonstrate the direct application of our polynomial Bernoulli number representation to $\zeta(3)$.
At first, we found
\begin{eqnarray}
\zeta (2n+1) =  (-)^{n+1} L_{2n} \frac{2(2\pi)^{2n}}{(2n)!}~,
\end{eqnarray}
and
\begin{eqnarray}
\zeta' (2n) = (-)^n \frac{(2\pi)^{2n}}{2(2n)!} B_{2n} \left( \frac{2n}{B_{2n}} L_{2n-1} -\gamma_E -ln(2\pi) \right)~, 
\end{eqnarray}
where $\gamma_E$ denotes the Euler constant. The L-numbers were first introduced by Bendersky in context with
the logarithmic Gamma-function \cite{Bendersky33,Kellner09}. The first derivative on the odd numbers $\zeta' (2n+1)$ may
be computed as follows
\begin{eqnarray}
\zeta'(3) = 2\pi^2 \zeta'' (-2) + \left(ln(2\pi) - \frac{3}{2} + \gamma_E \right) \zeta (3)~,
\end{eqnarray}
where
\begin{eqnarray}
\zeta''(-2) = 3L_2 - \frac{17}{108} + 4\sum_{n=2}^\infty \frac{B_{2n+2}H_{2n-2}}{(2n-1)(2n)(2n+1)(2n+2)}
\end{eqnarray}
can be written in terms of a Dirichlet series with the Bernoulli numbers and the harmonic numbers \cite{Spie90}
involved.  Furthermore, we recall one of our earlier results for $\zeta(3)$ published in \cite{Bentz94} where we
have introduced a fast converging series representation based on Benderskys L-numbers \cite{Bendersky33}
\begin{eqnarray}
\zeta(3)~=~\frac{\pi^2}{8}-\frac{\pi^2}{12}ln\left(\frac{\pi}{3}\right)+\frac{\pi^2}{3}\sum_{n=1}^\infty
\frac{\zeta(2n)}{2n(2n+1)(2n+2)}\left(\frac{1}{6}\right)^{2n}~.
\end{eqnarray}
Using the well-known Taylor-McLaurin series representation for $\zeta(s)$ \cite{Bor00}
\begin{eqnarray}
\zeta(s) = \frac{1}{s-1} + \frac{1}{2} + \sum_{n=2}^\infty \frac{B_k}{k} 
               \left(\begin{array}{c} 
                               s+k \\
                               k
               \end{array}\right),
\end{eqnarray}
with s$\in \mathbb -N$,
one finds for the first derivative at s=-2
\begin{eqnarray}
\zeta'(-2) =  -\frac{1}{36} +2\sum_{n=4}^\infty \frac{B_{2n}}{(2n-3)(2n-2)(2n-1)2n}~.
\end{eqnarray}
Following Bendersky \cite{Bendersky33}
the logarithmic gamma function ln$\Gamma_n$(x) of order n may be written as follows
\begin{eqnarray}
\lambda_n(x+1) = x^n ln(x) - L_n + ln\Gamma_n(x)~
\end{eqnarray}
where, for example, $\lambda_2(x+1)$ is defined as \cite{Bendersky33}
\begin{eqnarray}
\lambda_2(x+1)&=&\left( \frac{1}{3}x^3+\frac{1}{2}x^2+\frac{1}{6}x \right)~ln(x)
\nonumber \\ &-&\frac{1}{9}x^3+\frac{1}{12}x^2 + 2
\sum_{n=2}^\infty \frac{B_{2n+2}~x^{-2n+1}}{(2n-1)(2n)(2n+1)(2n+2)}~.
\end{eqnarray}
At x = 1 it follows
\begin{eqnarray}
\lambda_2(2)~ =~\zeta'(-2)~=~L_2~.
\end{eqnarray}
The L-numbers had been computed by Bendersky for all integer numbers n \cite{Bendersky33}. For n=2 we have
\begin{eqnarray}
L_2 = \frac{1}{48} \left(\frac{3}{2}-ln\left(\frac{\pi}{3}\right) + 4 \sum_{n=1}^\infty \frac{\zeta(2n)}{(2n)(2n+1)(2n+2)}~ 
\left(\frac{1}{6}\right)^{2n}\right)~.
\end{eqnarray}
This procedure allows a direct computation also of $\zeta'(-n)$ as shown by us before \cite{Bentz94}.
\begin{eqnarray}
\zeta' (-n) = - L_n + \frac{B_{n+1}}{n+1} \sum_{q=1}^n \frac{1}{q}~.
\end{eqnarray}
Within this procedure it is possible to find fast converging series in terms of the L-numbers for all $\zeta$(2n-1)
values. For example, for n=3 and n=5 one gets
\begin{eqnarray}
\zeta(5)&=&\frac{3\pi^2}{29}\zeta(3)-\frac{25\pi^4}{12528}+\frac{\pi^4}{1044}ln\left(\frac{\pi}{3}\right)
\nonumber \\ &-&\frac{4\pi^4}{87}\sum_{n=1}^\infty 
\frac{\zeta(2n)}{2n(2n+1)(2n+2)(2n+3)(2n+4)}\left(\frac{1}{6}\right)^{2n}~.
\end{eqnarray}
and
\begin{eqnarray}
\zeta(7)&=& \frac{72\pi^2}{659}\zeta(5)-\frac{2\pi^4}{1977}\zeta(3)+\frac{49\pi^6}{5337900}
+\frac{\pi^6}{266895}ln\left(\frac{\pi}{3}\right)
\nonumber \\ &-&\frac{32\pi^6}{5931}\sum_{n=1}^\infty 
\frac{\zeta(2n)}{2n(2n+1)(2n+2)(2n+3)(2n+4)(2n+5)(2n+6)}\left(\frac{1}{6}\right)^{2n}~.
\end{eqnarray}
For related constants like U$_2$ or U$_4$ we found:
\begin{eqnarray}
U_2&=&\sum_{n=1}^\infty (-)^{n+1}\frac{1}{(2n-1)^2} \nonumber \\
&=& 1~-~2 \sum_{n=1}^\infty  (-)^{n+1} \frac{nL_{2n}}{(2n)!}\left(\frac{\pi}{2}\right)^{2n}~,
\end{eqnarray}
and
\begin{eqnarray}
U_4&=&\sum_{n=1}^\infty (-)^{n+1}\frac{1}{(2n-1)^4} \nonumber \\
&=& \frac{9}{2}~-~2ln(2)~-~\frac{7}{6}U_2~-~\frac{7}{8}\zeta(3)~+~\frac{4}{3}\sum_{n=1}^\infty (-)^{n+1}
\frac{n^3L_{2n+2}}{(2n+2)!}\left(\frac{\pi}{2}\right)^{2n+2}~.
\end{eqnarray}
These examples show that not only zeta values can be expressed by Bendersky's L-numbers and as a consequence
by the even Bernoulli numbers. With this work we introduce a new polynomial representation of the even Bernoulli
numbers which allows for a fast computation of these constants, which is comparable with the computation by
use of available BBP-like formulas. Furthermore, our method allows for the use of polylogarithmic identities
on which BBP formulas are typically based.

The paper is organized as follows. In the next section we introduce a new formula for the even Bernoulli
numbers which serves as a basis for our polynomial representation. We present a prove of this formula
discuss it in context with other well-known series representations of Bernoulli numbers, giving a special
emphasize on the numerical effort which is needed to compute B$_{2n}$. In section 4 and 5 we present first
applications of our formalism to the computation of $\zeta$(3), $\zeta$(5) and $\zeta$(7) and compare with
corresponding BBP-type formulas concerning the numerical effort. The final section 6 gives a summary of the
results and a short outlook. In the appendix section we present as a further application of our polynomial
representation the computation of B$_{2n}$ as a function of B$_{2n-2}$ only. Based on this formula an even
faster computation of zeta numbers by one order of magnitude is possible.

\section{An alternative formula for the even Bernoulli numbers}
In this section we present the basic equation from which the Bernoulli numbers B$_{2n}$ can be computed. \\

\noindent
{\bf \large Theorem 2.1} \\ 

It follows for $B_{2n}$ 
\begin{eqnarray}
 B_{2n}&=&\frac{2^{2n+1}(2n)!}{2^{2n}-2}~\sum^{n}_{k=1} \frac{(-)^{k+1}}{(2n+k)!}
 \left(\begin{array}{c} 
            n+1 \\
            k+1
 \end{array}\right) \nonumber \\ &*&
~\sum^{[\frac{k}{2}]}_{l=0} (-)^l
 \left(\begin{array}{c} 
            k \\
            l
 \end{array}\right) (\frac{k}{2}-l)^{k+2n}~.
\end{eqnarray}
The formula shown above represents an alternative computational scheme for the Bernoulli numbers, as only the
even numbers are calculated, where the outer summation index runs up to n only. All other formulas available from
literature \cite{Comt74,Gould72,Dil89,Todo91} including the representation which was introduced in 2009 \cite{irak09}
compute the Bernoulli numbers B$_{n}$ by a summation up to n. A direct computation of the Bernoulli numbers by use
of Eq.~(2.1) is possible, where a runtime comparison with a variety of other existing recurrence formula reveals
that our formula  allows for the computation of B$_{2n}$ numbers which is typically faster by factors ranging
between 1.5 and 5.4. Corresponding numerical results are presented in Table I: Furthermore,
and this is the most important point, our double-sum like representation allows a direct link to the so called central factorial
numbers described in detail by Riordan \cite{Riordan68}, which finally offers the possibility to introduce a polynomial
representation of the Bernoulli numbers not known before. We show later on that this polynomial representation allows for a
fast computation of all odd zeta-numbers $\zeta(2n+1)$ and related constants discussed in the introduction.
 
To prove the theorem we first list some helpful identities.
\begin{table} 
\begin{center}
\begin{tabular}[bp]{||c||c|c|c|c|c|c||}
   \hline
   &&&&&&\\
              & \cite{x1} & \cite{x2} & \cite{x3} & \cite{Gould72} & \cite{x5} & This work  \\ 
   &&&&&&\\
   \hline
   B$_{10}$ & 2.5 & 1.5 & 1.5 & 3.0 & 2.0 & 1.0 \\ 
   \hline
   B$_{20}$ & 5.0 & 1.8 & 2.0 & 5.2 & 3.0 & 1.0 \\ 
   \hline
   B$_{30}$ & 4.8 & 1.6 & 1.7 & 5.3 & 3.1 & 1.0 \\ 
   \hline
   B$_{40}$ & 4.8 & 1.6 & 1.8 & 5.4 & 3.1 & 1.0 \\ 
   \hline
\end{tabular}
\caption{Runtime comparison between different recurrence relations for Bernoulli numbers with our present formula. Shown is
the factor which results from an explicit calculation of the even Bernoulli numbers for n=10, 20, 30 and 40.}
\end{center}
\end{table} 

\begin{eqnarray}
\frac{k-2l}{i+1}
\left(\begin{array}{c} 
            i+1 \\
            k-l
 \end{array}\right)
\left(\begin{array}{c} 
            k-l \\
            l
 \end{array}\right)~=~
\left(\begin{array}{c} 
            i \\
            k-l-1
 \end{array}\right)
\left(\begin{array}{c} 
            k-l-1 \\
            l
 \end{array}\right)
\end{eqnarray}~.

This identity is simply proved by a direct calculation of both sides.

\begin{eqnarray}
\sum_{i=0}^{n-1}
\left(\begin{array}{c} 
            k-2l \\
            i
 \end{array}\right)
\left(\begin{array}{c} 
            n \\
            i
 \end{array}\right)^{-1}=
\frac{\left( n+1-
\left(\begin{array}{c} 
            k-2l \\
            n
 \end{array}\right)
\right)}{n+1-k+2l}~.
\end{eqnarray}

Identity (2.3) is found in \cite{Riordan68,Gould721}.

\begin{eqnarray}
\frac{n+1}{n+1-k+2l}
\left(\begin{array}{c} 
            n \\
            k-l
 \end{array}\right)
\left(\begin{array}{c} 
            k-l \\
            l
 \end{array}\right) &=& \nonumber \\ &&
\left(\begin{array}{c} 
            n+1 \\
            k-2l
 \end{array}\right)
\left(\begin{array}{c} 
            n-k+2l \\
            l
 \end{array}\right)~.
\end{eqnarray}

Again, identity (2.4) is proved by a direct calculation of both sides. 
Now we start our proof with the well-known recursive formula \cite{Hu13}: 

\begin{eqnarray}
B_{2n}~=~\frac{1}{2}-\frac{1}{2n+1}\sum_{l=0}^{n-1}
\left(\begin{array}{c} 
            2n+1 \\
            2l
 \end{array}\right) B_{2l}~.
\end{eqnarray}
This recursive relation can be used to formulate an iterative solution for the Bernoulli numbers
\begin{eqnarray}
S_{i+1}(n)~=~\sum^{n-1}_{k=i}
\left(\begin{array}{c} 
            2n+i+1 \\
            2k+i
\end{array}\right)~S_i(k)~,
\end{eqnarray}
i=1, 2, 3, ...$\in \mathbb N$ where it follows for S$_1(n)$
\begin{eqnarray}
S_1(n)~=~1+\frac{1}{2}\sum_{l=1}^{n-1}
\left(\begin{array}{c} 
            2n+1 \\
            2l
\end{array}\right)~=~\frac{1}{2}(2^{2n}-2n)~.
\end{eqnarray}
In second order this recursive formula gives
\begin{eqnarray}
S_2(n)=~\frac{9}{8}3^{2n}-(2n+1)2^{2n}+2n^2+n-\frac{9}{8}~.
\end{eqnarray}
Then the Bernoulli numbers B$_{2n}$ result in 
\begin{eqnarray}
\frac{2B_{2n}}{2n}&=&\frac{1}{2n}-\frac{2S_1(n)}{2n(2n+1)}+\frac{2S_2(n)}{2n(2n+1)(2n+2)} \nonumber \\
&-&\frac{2S_3(n)}{2n(2n+1)(2n+2)(2n+3)}\pm ...
\end{eqnarray}
For a direct computation of (2.9) we convert in a first step all partial sums S$_i(n)$,~i=1, 2, 3,... $\in \mathbb N$
to a more appropriate form. For $i=1$, 2 and 3 it follows explicitly
\begin{eqnarray}
S_1(n)~=~\frac{1}{4}2^{2n+1}-\frac{1}{2}(2n)1^{2n+1}~,
\end{eqnarray}
\begin{eqnarray}
S_2(n)&=&\frac{1}{8}3^{2n+2}+\frac{1}{4}2^{2n+2}-\frac{1}{8}1^{2n+2}-2(2n+2)S_1(n)\nonumber \\ 
&-&\frac{1}{2}(2n+1)(2n+2)~,
\end{eqnarray}
\begin{eqnarray}
S_3(n)&=&\frac{1}{16}4^{2n+3}+\frac{1}{8}3^{2n+3}-\frac{1}{8}2^{2n+3}-\frac{3}{8}1^{2n+3}
-3(2n+3)S_2(n) \nonumber \\ &-& 3(2n+2)(2n+3)S_1(n)-\frac{1}{2}(2n+1)(2n+2)(2n+3)~.
\end{eqnarray}
Formula (2.12) represents a ternary sum because of the repeated use of the
recursion formula (2.6), and as a consequence S$_n(n)$ denotes a n-fold sum. At a first glance, this fact
sounds very discouraging but fortunately most of the terms cancel if all partial sums $S_i(n)$ with
index i=1,...,n-1 will be inserted in the corresponding n-fold partial sum S$_n(n)$ with index n. To
illustrate this procedure we compute the partial sum with index i=3 \\

\noindent
{\bf \large Example 2.1} \\
 
\begin{eqnarray}
S_3(n)&=&\frac{1}{16}4^{2n+3}+\frac{1}{8}3^{2n+3}-\frac{1}{8}2^{2n+3}-\frac{3}{8}1^{2n+3} \nonumber \\
&-&3(2n+3)\left(\frac{1}{8}3^{2n+2}+\frac{1}{4}2^{2n+2}-\frac{1}{8}1^{2n+2}\right) \nonumber \\
&+&3(2n+2)(2n+3)\left(\frac{1}{4}2^{2n+1}-\frac{1}{2}(2n)1^{2n+1}\right)  \nonumber \\
&-&\frac{1}{2}(2n+1)(2n+2)(2n+3)~.
\end{eqnarray}

This algorithm has to be applied to each partial sum S$_i(n)$ in a recursive sense. The result
for S$_i(n)$ is \\

\noindent
{\bf \large Lemma 2.1} \\

\begin{eqnarray}
S_i(n)&=&\frac{1}{2^{i+1}}\sum^i_{k=0}(-)^i(i+1-k)^{2n+i}\sum^k_{l=0}h(l)~(k-l)! \nonumber \\
&*& \left(\frac{2}{i+1-k}\right)^{k-l}
\left(\begin{array}{c} 
            2n+i \\
            k-l
\end{array}\right)
\left(\begin{array}{c} 
            i \\
            i+1-k
\end{array}\right)
\end{eqnarray}

with
\begin{eqnarray}
h(l)&=& \left(cos\frac{\pi l}{2}-sin\frac{\pi l}{2} \right)
\left(\begin{array}{c} 
            i+l-k \\
            t_1(l) 
\end{array}\right) 
\nonumber \\ &-& \left(cos\frac{\pi l}{2}+sin\frac{\pi l}{2}\right)
\left(\begin{array}{c} 
            i+l-k \\
           t_2(l)
\end{array} \right)~,
\end{eqnarray}
with $t_1(l)=\frac{l}{2}-\frac{1}{4}+\frac{(-)^l}{4})$ and $t_2(l)=\frac{l}{2}-\frac{3}{4}-\frac{(-)^l}{4})$.
The index t$_1(l)$ produces the numbers 0, 0, 1, 1, 2, 2, 3, 3, for $l$=0, 1, 2, 3, ..., k. Thus has the same
values for l even or odd. The second index t$_2(l)$ is slightly different as this sequences starts with -1 for
l=0. The proof can done by induction. 

As an example, we show in table II the complete list of numbers which result for $S_i(n)$,
where i=1,...,6 and n=1,...,10:
\begin{table} 
\begin{center}
\begin{tabular}[bp]{||c|c|c|c|c|c|c|c||}
   \hline
   &&&&&&\\
           & S$_1$(n) & S$_2$(n) & S$_3$(n) & S$_4$(n) & S$_5$(n) & S$_6$(n) \\ 
   &&&&&&\\
   \hline
   n=1  & 1 & 0 & 0 & 0 & 0 & 0  \\ 
   \hline
   n=2  & 6 & 20 & 0 & 0 & 0 & 0 \\ 
   \hline
   n=3  & 29 & 392 & 1680 & 0 & 0 & 0 \\ 
   \hline
   n=4  & 124 & 5112 & 73920 & 369600 & 0 & 0 \\ 
   \hline
   n=5  & 507 & 55220 & 2000856 & 30270240 & 16816800 & 0 \\ 
   \hline
   n=6  & 2024 & 544700 & 43099056 & 1462581120 & 2.287$^{+10}$ & 1.372$^{+11}$ \\ 
   \hline
   n=7  & 8185 & 5135184 & 821292576 & 5.475$^{+10}$ & 1.788$^{+12}$ & 2.868$^{+13}$ \\ 
   \hline
   n=8  & 32760 & 47313584 & 1.459$^{+10}$ & 1.772$^{+12}$ & 1.062$^{+14}$ & 3.384$^{+15}$ \\ 
   \hline
   n=9  & 131063 & 430867484 & 2.489$^{+11}$ & 5.250$^{+13}$ & 5.365$^{+15}$ & 2.991$^{+17}$ \\ 
   \hline
   n=10 & 524278 & 3900612564 & 4.140$^{+12}$ & 1.470$^{+15}$ & 2.444$^{+17}$ & 2.218$^{+19}$ \\ 
   \hline
\end{tabular}
\caption {Computed numbers for $S_i(n)$ with i running from 1 to 6 and n running from 1 to 10.}
\end{center}
\end{table}

\noindent
{\it Proof.} \\

Using now (2.9), doing few elementary manipulations, one ends up with a
first expression for B$_{2n}$ in form of a ternary sum only
\begin{eqnarray}
B_{2n}&=&2^{2n}(2n)!~\sum_{i=0}^n\sum_{k=0}^i\sum_{l=0}^{t_1(k)}(-)^{k+i+l}~\frac{2n+i-2k+4l+1}{i-k+l+1}
\nonumber \\ &*& 
\left(\begin{array}{c} 
            i \\
            k-2l
 \end{array}\right)
\left(\begin{array}{c} 
            i-k+2l \\
            l
 \end{array}\right)
\frac{\left( \frac{i+1-k}{2} \right)^{2n+i-k+2l+1}}{(2n+i-k+2l+1)!}~,
\end{eqnarray}
Making use of lemma (2.1) the different terms in (2.16) can be rearranged, leading to a simpler formula for $B_{2n}$
\begin{eqnarray}
B_{2n}&=& \frac{2^{2n+1}(2n)!}{2^{2n}-2}\sum_{i=1}^n\sum_{k=0}^i\sum_{l=0}^{t_1(k)}(-)^{k+i+l}
\nonumber \\ &*&
\left(\begin{array}{c} 
            i \\
            k-l+1
 \end{array}\right)
\left(\begin{array}{c} 
            k-l-1 \\
            l
 \end{array}\right)
\frac{\left( \frac{i+1-k}{2} \right)^{2n+i-k+2l+1}}{(2n+i-k+2l+1)!}~.
\end{eqnarray}
A further rearrangement results in
\begin{eqnarray}
B_{2n}&=& -\frac{2^{2n}(2n)!}{2^{2n-1}-1}\sum_{k=0}^{n-1}\sum_{l=0}^{t_1(k)}(-)^{n-k+l}
\nonumber \\ &*&
\left(\begin{array}{c} 
            n \\
            k-l
 \end{array}\right)
\left(\begin{array}{c} 
            k-l \\
            l
 \end{array}\right)
\frac{\left( \frac{n-k}{2} \right)^{3n-k+2l}}{(3n-k+2l)!}~
\left[ \sum_{i=0}^{n-1}
\left(\begin{array}{c} 
            k-2l \\
            i
 \end{array}\right)
\left(\begin{array}{c} 
            n \\
            i
 \end{array}\right)^{-1} \right]~,
\end{eqnarray}
where we have interchanged the sum over i with the other two sums. As the variable i appears in the
inner sum only this offers the possibility to reduce the ternary sum to a binary one. This is done
by use of lemma(2.2) \cite{Riordan68,Gould721}.
 
As the contribution from the Binomial coefficient 
$\left(\begin{array}{c} 
            k-2l \\
            n
 \end{array}\right)$ is zero it remains for $B_{2n}$ 
\begin{eqnarray}
B_{2n}&=& (-)^{n+1}\frac{2^{2n}(2n)!}{2^{2n-1}-1}\sum_{k=0}^{n-1}\sum_{l=0}^{t_1(k)}(-)^{l-k}
\nonumber \\ &*&
\left(\begin{array}{c} 
            n+1 \\
            k-2l
 \end{array}\right)
\left(\begin{array}{c} 
            n-k+2l \\
            l
 \end{array}\right)
\frac{\left( \frac{n-k}{2} \right)^{3n-k+2l}}{(3n-k+2l)!}~,
\end{eqnarray}
where we made use of lemma (2.3). Next, we rewrite (2.19) as follows
\begin{eqnarray}
B_{2n}&=&(-)^{n+1}\frac{2^{2n}(2n)!}{2^{2n-1}-1}\sum_{k=0}^n \frac{(-)^k}{(3n-k)!}
\left(\begin{array}{c} 
            n+1 \\
            k
 \end{array}\right) \nonumber \\ &*&
\sum_{l=0}^{t_1(n-k)} (-)^l
\left(\begin{array}{c} 
            n-k \\
            l
 \end{array}\right)
\left( \frac{n-k}{2}-l \right)^{3n-k}~,
\end{eqnarray}
and use the relation

\begin{eqnarray}
\sum_{l=0}^{t_1(n-k)} (-)^l
\left(\begin{array}{c} 
            n-k \\
            l
 \end{array}\right)
\left( n-k-2l \right)^{3n-k} &=& \nonumber \\ &&
\sum_{l=0}^{t_1(k)} (-)^l
\left(\begin{array}{c} 
            k \\
            l
 \end{array}\right)
\left( k-2l \right)^{2n+k}~,
\end{eqnarray}

which is simply proved in replacing n-k by k as the summation is symmetric in these indices. Finally we result in
\begin{eqnarray}
B_{2n}&=&\frac{2^{2n}(2n)!}{2^{2n-1}-1}\sum_{k=0}^n(-)^{k+1}
\left(\begin{array}{c} 
            n+1 \\
            n-k
 \end{array}\right)
\frac{(\frac{1}{2})^{2n+k}}{(2n+k)!} \nonumber \\ &*&
\sum_{l=0}^{t_1(k)} (-)^l
\left(\begin{array}{c} 
            k \\
            l
 \end{array}\right)
\left( k-2l \right)^{2n+k}~.
\end{eqnarray}
Finally the summation over $t_1(k)$ can be replaced by a summation over $[\frac{k}{2}]$ because the inner sum behaves
symmetric with respect to the index $t_1(k)$. This last step completes the proof of the theorem as the binomial coefficients
in (2.1) and (2.22) are identical. $\Box$ 

\section{The $\alpha$-polynomial as a generator to compute Bernoulli numbers}
We define the inner sum in (2.1) as a polynomial of order n in the variable k $\in \mathbb N$ with coefficients
$\alpha^{(n)}_l$, l=1,...,n, which can be directly computed in a numerical sense by use of the definition shown
below. To define the inner sum as a polynomial at a first glance seems to be quite unconventional, but this 
definition is based on corresponding relations which can be found, for example, in \cite{Riordan68,Gould721} for
k=0,1 and 2. The generalization to k$\ge$3 is straightforward, where explicit formulas for the coefficients
$\alpha^{(n)}_l$ will be given at the end of this section. We have then \\

\noindent
{\bf \large Proposition 3.1} \\

\begin{eqnarray}
\sum^{n}_{l=1} \alpha^{(n)}_l~k^l~=~\frac{2(2n)!}{(2n+k)!}~\sum^{[\frac{k}{2}]}_{j=0} (-)^j  
 \left(\begin{array}{c} 
            k \\
            j
 \end{array}\right) (\frac{k}{2}-j)^{k+2n}~.
\end{eqnarray}

Therefore, it follows for $B_{2n}$
\begin{eqnarray}
 B_{2n} = \frac{2^{2n}}{2^{2n}-2}~ \sum^{n}_{k=1} (-)^{k+1}
 \left(\begin{array}{c} 
            n+1 \\
            k+1
 \end{array}\right) \sum^{n}_{l=1} \alpha^{(n)}_l~k^l~.
\end{eqnarray}
Next, we show that the $\alpha$-polynomial can be defined for negative integer numbers. Furthermore,
it will be demonstrated that a computation of the $\alpha$-polynomial for negative integers provides
a new option to define the Bernoulli numbers. We start with a recursive formula which the $\alpha$-polynomial
fulfills. It is \\

\noindent
{\bf \large Lemma 3.1} \\

\begin{eqnarray}
\sum_{l=1}^n \alpha_l^{(n)} k^l &=& \frac{k(k-1)}{(2n+k)(2n+k-1)} \sum_{l=1}^n \alpha_l^{(n)} (k-2)^l
\nonumber \\ &+& \frac{1}{4}k^2 \frac{2n(2n-1)}{(2n+k)(2n+k-1)} \sum_{l=0}^{n-1} \alpha_l^{(n-1)} k^l
\end{eqnarray}

with $\alpha^{0}_0 = 1.$ 

To prove the above lemma let us first recall the following recursive formula obtained by Riordan
\cite{Riordan68}, for the so called central factorial numbers \footnote{Riordan discusses these numbers
in context with the so called central difference operator $\delta$ in chapter V of his book 'Combinatorial
identities'.}

\begin{eqnarray}
T(n,k)~=~T(n-2,k-2)~+~\frac{1}{4}k^2T(n-2,k)~,
\end{eqnarray}
with
\begin{eqnarray}
T(n,k)~=~\frac{1}{k!}\sum_{l=0}^k (-)^l
\left(\begin{array}{c} 
            k \\
            l
 \end{array}\right)
\left(\frac{k}{2}-l \right)^n~. 
\end{eqnarray}
From this it follows immediately
\begin{eqnarray}
T(2n+k,k)~=~
\left(\begin{array}{c} 
            2n+k \\
            k
 \end{array}\right)
\sum_{l=1}^n \alpha_l^{(n)}k^l~.
\end{eqnarray}
Now it is easy to complete the proof. First we write

\begin{eqnarray}
T(2n+k,k)~=~T(2n+k-2,k-2)~+~\frac{1}{4}k^2T(2n+k-2,k)~.
\end{eqnarray}
Substituting here (3.6) the recursive formula for the $\alpha$-polynomial follows directly for k$\in \mathbb N$.
Furthermore, we have shown that the central factorial numbers are closely related to the $\alpha$-polynomial,
which can be regarded as a key quantity in computing the Bernoulli numbers.

Using the notation \\

\noindent
{\bf \large Definition 3.1} \\

\begin{eqnarray}
A^{(n)}(k) = \sum_{l=1}^n \alpha_l^{(n)} k^l~,
\end{eqnarray}

it follows for example for k=1 \\

\noindent
{\bf \large Example 3.1} \\

\begin{eqnarray}
A^{(n)}(1) = \frac{1}{2^{2n}~(2n+1)}.
\end{eqnarray}

This way the sum $A^{(n)}(k)$ may be calculated for all k $\in \mathbb N$. To extend the calculational
scheme to negative integers we introduce the following relation \\

\noindent
{\bf \large Lemma 3.2} \\

\begin{eqnarray}
A^{(n)}(-k)~=~-k
               \left(\begin{array}{c} 
                               n + k \\
                               n
               \end{array}\right)
\sum_{l=1}^n \frac{(-)^{l+1}}{l+k} 
               \left(\begin{array}{c} 
                               n \\
                               l
               \end{array}\right)
A^{(n)}(l)~.
\end{eqnarray}

To prove the above relation we use the following identity \cite{Riordan68,Gould721} \\

\noindent
{\bf \large Proposition 3.2} \\

\begin{eqnarray}
\sum_{l=1}^n \frac{(-)^l}{k+l}
\left(\begin{array}{c} 
                               n \\
                               l
               \end{array}\right)
l^m~=~(-)^m~k^{m-1}
\left(\begin{array}{c} 
                               n+k \\
                               k
               \end{array}\right)^{-1}~.
\end{eqnarray} \\

\noindent
{\it Proof.} \\

Writing now
\begin{eqnarray}
\sum_{m=1}^n \alpha_m^{(n)}(-k)^m&=& -k
\left(\begin{array}{c} 
                               n+k \\
                               k
               \end{array}\right)
\sum_{l=1}^n \frac{(-)^{l+1}}{l+k}\sum_{m=1}^n \alpha_m^{(n)}l^m \nonumber \\ &=&
k 
\left(\begin{array}{c} 
                               n+k \\
                               k
               \end{array}\right)
\sum_{m=1}^n \alpha_m^{(n)} \sum_{l=1}^n \frac{(-)^l}{k+l}
\left(\begin{array}{c} 
                               n \\
                               l
               \end{array}\right) l^m \nonumber \\ &=&
\sum_{m=1}^n \alpha_m^{(n)}(-k)^m~.
\end{eqnarray}~~$\Box$

As a consequence we find the result that $A^{(n)}(-1)$ is "proportional" to B$_{2n}$ \\

\noindent
{\bf \large Example 3.2} \\

\begin{eqnarray}
A^{(n)}(-1) = -\frac{2^{2n}-2}{2^{2n}}~B_{2n}~.
\end{eqnarray}
Using (3.3) it follows further
\begin{eqnarray}
A^{(n)}(-2)~=~(2n-1)~B_{2n}~, 
\end{eqnarray} 
and for example
\begin{eqnarray}
A^{(n)}(-3)~&=&~\frac{1}{2}(2n)(2n-1)\frac{2^{2n}-2}{2^{2n}}~B_{2n} \nonumber \\
&-&\frac{1}{4}(2n-1)(2n-2)\frac{2^{2n-2}-2}{2^{2n-2}}~B_{2n-2}~, 
\end{eqnarray} 

where in the expression for $A^{(n)}(-3)$ both B$_{2n}$ and B$_{2n-2}$ appear. More detailed information
about the properties of the A$^{(n)}(k)$ can be obtained by changing to an alternative representation, which
was first introduced by Riordan \cite{Riordan68}, again in context with the so called central factorial numbers \\

\noindent
{\bf \large Lemma 3.3} \\
\begin{eqnarray}
A^{(n)}(k) = \frac{1}{2^{2n}} \left(\begin{array}{c}
                                                      2n + k \\
                                                      2n
                                                  \end{array}\right)^{-1}
\sum_{l=0}^n a_{n,l} \left(\begin{array}{c}
                               2n + k \\
                               2n + l
                    \end{array}\right)~,
\end{eqnarray}
with a$_{n,0}=\delta_{n,0}$ and a$_{n,1}$=1~. \\

\noindent
{\it Proof.} \\

Following Riordan \cite{Riordan68} we first have
\begin{eqnarray}
2^{2n} T(2n+k,k)~=~\sum_{l=0}^n a_{n,l}
\left(\begin{array}{c}
                               2n + k \\
                               2n + l
                    \end{array}\right)~.
\end{eqnarray}
Using (3.6) it follows
\begin{eqnarray}
T(2n+k,k)&=&\frac{1}{2^{2n}} \sum_{l=0}^k a_{n,l}
\left(\begin{array}{c}
                               2n + k \\
                               2n + l
                    \end{array}\right) \nonumber \\ &=&
\left(\begin{array}{c}
                               2n + k \\
                               2n 
                    \end{array}\right) A^{(n)}(k)
\end{eqnarray}
and this completes the proof.~$\Box$

The coefficients $a_{n,l}$ are defined by the recursive formula presented below \cite{Riordan68}
\begin{eqnarray}
a_{n,l} = \sum^{n-1}_{i=l-1} a_{i,l-1}
               \left(\begin{array}{c} 
                               2n + l - 1 \\
                               2n -2i
               \end{array}\right)~,
\end{eqnarray}
with a$_{n,1}$=a$_{0,0}$ = 1 and a$_{n,0}$ = 0 for n$>$0.
It should be mentioned here that no obvious connection between the central factorial numbers and
the Bernoulli numbers exists. This is because the recursive relation (3.19), which allows to compute
the central factorial numbers is linked by (3.16) only to the $\alpha$-polynomials but not directly
to the Bernoulli numbers. Thus, it becomes clear that the Bernoulli numbers were not discussed by
Riordan in context with the central factorial numbers. As an example it follows for n=2 \\

\noindent
{\bf \large Example 3.3} \\
 
\begin{eqnarray}
A^{(2)}(k) &=& \alpha_1^{(2)}k + \alpha_2^{(2)}k^2 = \frac{1}{2^4} \left(\begin{array}{c}
                                                                     k + 4 \\
                                                                     k
                                                                  \end{array}\right)^{-1}
\left[\left(\begin{array}{c}
                    k + 4 \\
                    5
       \end{array}\right)
+ 10
      \left(\begin{array}{c}
                    k + 4 \\
                    6
       \end{array}\right)
\right] \nonumber \\ &=& \frac{1}{48} k^2 + \frac{1}{4} B_4 k~,
\end{eqnarray}
where the Bernoulli number B$_4$ appears as the coefficient of the lowest order polynomial term. 

Computing (3.16) for k = -1 we find for B$_{2n}$
\begin{eqnarray}
B_{2n}~=~\frac{1}{2^{2n}-2}~\sum_{l=1}^n~(-)^{l+1}
~a_{n,l}~
               \left(\begin{array}{c} 
                               2n+l \\
                               2n
               \end{array}\right)^{-1}~. 
\end{eqnarray}
In a next step we substitute for $a_{n,l}$ \\

\noindent
{\bf \large Definition 3.2} \\

\begin{eqnarray}
a_{n,l}~=~\frac{(2n)!}{6^n}l
               \left(\begin{array}{c} 
                               2n+l \\
                               2n
               \end{array}\right)P^{(n+1-l)}(n)
\end{eqnarray}

by introducing a new type of polynomials. The result for B$_{2n}$ is then
\begin{eqnarray}
B_{2n}~=~\frac{(2n)!}{(2^{2n}-2)6^n}~\sum_{l=1}^n~(-)^{l+1}~l~P^{(n+1-l)}(n)~.
\end{eqnarray}
The P-polynomials can be obtained from (3.22) but they are also connected to the
$\alpha$-polynomials \\

\noindent
{\bf \large Lemma 3.4} \\

\begin{eqnarray}
P^{(n+1-l)}(n)~=~\frac{2^{2n}6^n(-)^l}{l(2n)!}~\sum^l_{k=1}~(-)^k
               \left(\begin{array}{c} 
                               l \\
                               k
               \end{array}\right) 
A^{(n)}(k)~,
\end{eqnarray}
and vice versa for $A^{(n)}(k)$
\begin{eqnarray}
A^{(n)}(k)~=~\frac{(2n)!}{2^{2n}6^n}~\sum^n_{l=1}~l
               \left(\begin{array}{c} 
                               k \\
                               l
               \end{array}\right)
P^{(n+1-l)}(n)~.
\end{eqnarray} \\

\noindent
{\it Proof.} \\

First, the relation for $A^{(n)}(k)$ results by substituting (3.22) in (3.16) and by reformulating
the product of the three binomial coefficients. The above formula for the P-polynomials is then
proved by substituting (3.24) in (3.23). This gives 
\begin{eqnarray}
B_{2n}~=~\frac{2^{2n}}{2^{2n}-2}\sum_{l=1}^n\sum_{k=1}^l (-)^{k+1}
\left(\begin{array}{c} 
                               l \\
                               k
               \end{array}\right)
A^{(n)}(k)~.
\end{eqnarray}  
Replacing now $l$ by $n+l-1$ it follows 
\begin{eqnarray}
B_{2n}~=~\frac{2^{2n}}{2^{2n}-2}\sum_{l=1}^n\sum_{k=1}^{n+1-l} (-)^{k+1}
\left(\begin{array}{c} 
                               n+1-l \\
                               k
               \end{array}\right)
A^{(n)}(k)~.
\end{eqnarray}  
The double sum can be written as
\begin{eqnarray}
\sum_{k=1}^n (-)^{k+1}
\left(\begin{array}{c} 
                               n \\
                               k
               \end{array}\right)
A^{(n)}(k)+...+
\sum_{k=1}^n (-)^{k+1}
\left(\begin{array}{c} 
                               1 \\
                               k
               \end{array}\right)
A^{(n)}(k) \nonumber \\
= \sum_{k=1}^n (-)^{k+1}
\left(\begin{array}{c} 
                               n+1 \\
                               k+1
               \end{array}\right)
A^{(n)}(k)~,
\end{eqnarray}
which finally gives (3.2).~$\Box$

As an example, we have for the first three P-polynomials \\

\noindent
{\bf \large Example 3.4} \\

\begin{eqnarray}
P^{(1)}(n)~=~\frac{1}{n},~~~~P^{(2)}(n)~=~\frac{3}{10},~~~~P^{(3)}(n)~=~\frac{3(21n-43)}{2^3~5^2~7}~.
\end{eqnarray}
All P-polynomials can be regarded as ordinary polynomials, apart from the first one.

It was shown by (3.10) that the $\alpha$-polynomials can be defined for negative integers.
The corresponding formula computed as a function of the P-polynomials is given by \\

\noindent
{\bf \large Lemma 3.5} \\

\begin{eqnarray}
A^{(n)}(-k)~=~(-)^n \frac{(2n)!}{6^n2^{2n}}~k~\sum_{l=1}^n (-)^{l+1}
               \left(\begin{array}{c} 
                               n+k-l\\
                               k
               \end{array}\right)~P^{(l)}(n)~.
\end{eqnarray} \\

\noindent
{\it Proof.} \\

Computing (3.25) for k = -1 gives 
\begin{eqnarray}
A^{(n)}(-k)~=~\frac{(2n)!}{6^n2^{2n}}\sum_{l=1}^n (-)^l~l~
\left(\begin{array}{c} 
                               l+k-1\\
                               l
               \end{array}\right)~P^{(n+1-l)}(n)~.
\end{eqnarray}
Replacing again $l$ by $n+l-1$ it follows
\begin{eqnarray}
A^{(n)}(-k)&=&\frac{(2n)!}{6^n2^{2n}}\sum_{l=1}^n (-)^{n+1-l}(n+1-l)
	\left(\begin{array}{c} 
			n+k-l\\
			n+1-l
			\end{array}\right)~P^{(l)}(n) \nonumber \\ &=&
(-)^n\frac{(2n)!}{6^n2^{2n}}\sum_{l=1}^n (-)^{l+1}(n+k-l)
	\left(\begin{array}{c} 
			n+k-l-1\\
			k-1
			\end{array}\right)~P^{(l)}(n)~.
\end{eqnarray}~$\Box$ 

The relation (3.32) allows an explicit definition of the $\alpha^{(n)}_l$-coefficients
by comparing directly the corresponding coefficients in the polynomials on both sides. For
$\alpha^{(n)}_1$ we found the interesting result \\

\noindent
{\bf \large Corollary 3.1} \\

\begin{eqnarray}
\frac{B_{2n}}{2n}&=&\frac{(2n)!}{6^n 2^{2n}}~\sum_{l=1}^{n}(-)^{l+1}~
\frac{s(l,1)}{(l-1)!}~P^{(n+1-l)}(n)
\nonumber \\ &=&
\frac{(2n)!}{6^n 2^{2n}}~\sum_{l=1}^{n} (-)^{l+1}~P^{(n+1-l)}(n)~=~\alpha^{(n)}_1.
\end{eqnarray}
The coefficients for l$\ge$2 become 
\begin{eqnarray}
\alpha^{(n)}_2~=~\frac{(2n)!}{6^n 2^{2n}}~\sum_{l=2}^{n}
(-)^{l+2}~\frac{s(l,2)}{(l-1)!}~P^{(n+1-l)}(n)~,
\end{eqnarray}
\begin{eqnarray}
\alpha^{(n)}_n&=&\frac{(2n)!}{6^n 2^{2n}}~\sum_{l=n}^{n} (-)^{l+n}~
\frac{s(l,n)}{(l-1)!}~P^{(n+1-l)}(n) \nonumber \\ &=&
\frac{(2n)!}{6^n 2^{2n}}\frac{s(n,n)}{(n-1)!}P^{(1)}(n)~=~\frac{(2n)!}{6^n 2^{2n}n!}~,
\end{eqnarray}
where s(l,k) denote the Stirling numbers of first kind \cite{Niel1904,Zhao08,Mezo10}.
Formula (3.32) can be computed for different k-values. As an example we present here the computation
for k=0. The corresponding calculation for k=1 and 2 is shown in the appendix section. 
We find for k=0: \\

\noindent
{\bf \large Example 3.5} \\

\begin{eqnarray}
\sum^n_{l=1}(-)^{l+1}~P^{(n+1-l)}(n)~=~\frac{2^{2n}6^n}{2n(2n!)}B_{2n},
\end{eqnarray}
or
\begin{eqnarray}
\sum^n_{l=1}(-)^{l+1}~P^{(l)}(n)~=~(-)^{n+1}\frac{2^{2n}6^n}{2n(2n!)}B_{2n}~,
\end{eqnarray}
with
\begin{eqnarray}
\sum^n_{l=1}(-)^{l+1}~P^{(n+1-l)}(n)~=~(-)^{n+1}\sum^n_{l=1}(-)^{l+1}~P^{(l)}(n)~.
\end{eqnarray}

\section{Application to a series representation of $\zeta(3)$}
As mentioned in the introduction our polynomial representation of the Bernoulli numbers can be used to
compute rather fast converging sequences, for example, for $\zeta(3)$. Using (1.1) and (3.39) it follows 

\begin{eqnarray}
\zeta(2n)~=~\frac{1}{2}(2n)\zeta(2)^n\sum_{l=1}^n~(-)^{l+1}P^{(l)}(n).
\end{eqnarray}

In other words this special P-polynomial sum converges rather fast as n increases
\begin{eqnarray}
\lim_{n\rightarrow \infty}~\sum_{l=1}^n~(-)^{l+1}P^{(l)}(n)~=~0~.
\end{eqnarray}
This is essential because the prefactor in (3.39) growths much faster with n as the Bernoulli
numbers do. As a consequence the P-polynomial sum must compensate this behavior so that the
product results in $B_{2n}$. Also, a representation of $\zeta$(2n) in terms of $\zeta$(2) is possible
with the help of the polynomial representation. Furthermore, by use of (5.1) a direct representation
of the ln~sin(x) function is possible. We find
\begin{eqnarray}
ln~sin(\pi x)~&=&~ln(\pi x)~-~2\sum_{n=1}^{\infty} \frac{\zeta(2n)}{2n} \nonumber \\ 
&=& ln(\pi x)~-~\sum_{n=1}^{\infty}\zeta(2)^n \sum_{l=1}^n~(-)^{l+1}P^{(l)}(n)~x^{2n}~.
\end{eqnarray}
Following Bendersky \cite{Bendersky33} we have
\begin{eqnarray}
\zeta(3)~=~-6\pi^2~\int_0^{\frac{1}{6}}\int_0^{x}ln~sin(\pi y) dy -\frac{\pi^2}{12}~ln(2)~.
\end{eqnarray}
By integrating (5.3) twice and combining the result with (5.4) it follows
\begin{eqnarray}
\zeta(3)~=~\frac{\pi^2}{8}~-~\frac{\pi^2}{12}~ln\left(\frac{\pi}{3}\right)&+&\frac{\pi^2}{3} \sum_{n=1}^{\infty}
\frac{1}{2n(2n+1)(2n+2)}\left(\frac{\pi}{6\sqrt6}\right)^{2n} \nonumber \\ &-&
\frac{\pi^2}{20} \sum_{n=2}^{\infty} \frac{1}{(2n+1)(2n+2)}\left(\frac{\pi}{6\sqrt6}\right)^{2n}~+~R\left(\frac{\pi}{6\sqrt6}\right)~.
\end{eqnarray}
If we consider the first two sums only $\zeta(3)$ follows with an error $\delta <$ 10$^{-7}$. The convergence
is very fast, as higher order sums according to the polynomials P$^{(i)}$(n) with index i = 3,4,... start with
a summation index n = i instead of n = 1. A further advantage of this type of summation is that all sums appearing
can be expressed in terms of elementary functions based on logarithmic expressions. This allows a more
detailed insight on $\zeta(3)$ and in consequence on all other zeta-values evaluated at odd integer numbers.
In the next section we present an alternative series for  $\zeta(3)$, $\zeta(5)$ and $\zeta(7)$ which converge
even faster using Eq. (3.33) for k=2.

\section{Computation of $\zeta(3)$, $\zeta(5)$ and $\zeta(7)$} 
As the main application we use now an alternative polynomial representation to compute $\zeta(3)$, $\zeta(5)$ and $\zeta(7)$.
With Eq. (3.33) it follows: \\

\noindent
{\bf \large Proposition 4.1} \\
\begin{eqnarray}
ln~sin(\pi x)&=&ln(\pi x) \nonumber \\ &-& 
\sum_{n=1}^{\infty}\frac{2\zeta(2)^n}{(2n-1)2n} 
\left(\sum_{l=1}^n~(-)^{l+1}\left(\begin{array}{c} 
                               n+2-l\\
                               2
                        \end{array}\right)P^{(l)}(n) \right)~x^{2n}~.
\end{eqnarray}
As a consequence we find for $\zeta(3)$ by integration
\begin{eqnarray}
\zeta(3)~=~\frac{\pi^2}{8}~-~\frac{\pi^2}{12}~ln\left(\frac{\pi}{3}\right)~+~36\sum_{n=1}^{\infty}(-)^{n+1}c_n
\left(\frac{\pi}{6\sqrt6}\right)^{2n}~,
\end{eqnarray}
with
\begin{eqnarray}
c_1~=~\sum_{n=1}^{\infty}\frac{n(n+1)P^{(1)}(n)}{(2n-1)2n(2n+1)(2n+2)}\left(\frac{\pi}{6\sqrt6}\right)^{2n}~,
\end{eqnarray}
\begin{eqnarray}
c_2~=~\sum_{n=1}^{\infty}\frac{n(n+1)P^{(2)}(n+1)}{(2n+1)(2n+2)(2n+3)(2n+4)}\left(\frac{\pi}{6\sqrt6}\right)^{2n}~,
\end{eqnarray}
\begin{eqnarray}
c_3~=~\sum_{n=1}^{\infty}\frac{n(n+1)P^{(3)}(n+2)}{(2n+3)(2n+4)(2n+5)(2n+6)}\left(\frac{\pi}{6\sqrt6}\right)^{2n}~,
\end{eqnarray}
\begin{eqnarray}
c_4~=~\sum_{n=1}^{\infty}\frac{n(n+1)P^{(4)}(n+3)}{(2n+5)(2n+6)(2n+7)(2n+8)}\left(\frac{\pi}{6\sqrt6}\right)^{2n}~,
\end{eqnarray}
and this way for higher coefficients c$_n$.
Using these first four coefficients only we find the following approximate value for $\zeta(3)$:
\begin{eqnarray}
\zeta(3)~>~1.2020569031595738...~,
\end{eqnarray}
with an error $\delta <$ 0.2*10$^{-13}$. The convergence is about 3 orders of magnitude with the summation index n. This is
about 1 order of magnitude faster compared to the corresponding numerical values which result from the series representation
(1.7). Furthermore, this new series shows up with an alternating sign which has some benefit in estimating its the convergence
properties. A further advantage of this type of summation is that all infinite sums appearing can be expressed in terms of
elementary functions based on logarithmic expressions. This allows a more detailed insight on $\zeta(3)$ and in consequence
on all other zeta-values evaluated at odd integer numbers. The explicit comparison is shown in Table III:
\begin{table}
\begin{center}
\begin{tabular}[t]{||c||c|c|c||}
   \hline
   &&&\\
   $\zeta$(3) & $\zeta$(3)-Zeta series & $\zeta$(3)-(Zeta series+BP) & $\zeta$(3)-BBP formula \cite{kun10} \\ 
   &&&\\
   \hline
   1st order &&& \\
   n=1 & $\delta$=0.2*10$^{-04}$ &  $\delta$=0.2*10$^{-05}$ & $\delta$=0.5*10$^{-06}$ \\ 
   \hline
   2nd order &&& \\
   n=2 & $\delta$=0.2*10$^{-06}$ &  $\delta$=0.1*10$^{-08}$ & $\delta$=0.9*10$^{-10}$ \\ 
   \hline
   3rd order &&& \\
   n=3 & $\delta$=0.3*10$^{-08}$ &  $\delta$=0.8*10$^{-12}$ & $\delta$=0.4*10$^{-13}$ \\ 
   \hline
   4th order &&& \\
   n=4 & $\delta$=0.4*10$^{-10}$ &  $\delta$=0.2*10$^{-13}$ & $\delta$=0.2*10$^{-16}$ \\ 
   \hline
\end{tabular}
\caption{Approximate computation of $\zeta$(3) as a function of the summation index n by
use of Eq.~(1.7) without and with use of the polynomial representation. The numerical errors are compared
to the BBP-type formula \cite{kun10}}
\end{center}
\end{table}
where the numerical comparison has been performed between our approach and the BBP-formalism by use of the following
formula \cite{kun10}:
\begin{eqnarray}
\frac{5\pi^2ln(3)}{104}-\frac{3ln(3)^3}{104} &=& \frac{1}{1053}\sum_{n=1}^{\infty} \left( \frac{1}{729} \right)^n
\Big( \frac{729}{(12n+1)^3}+\frac{243}{(12n+2)^3}-\frac{81}{(12n+4)^3} \nonumber \\ &-&
\frac{81}{(12n+5)^3}-\frac{54}{(12n+6)^3}-\frac{27}{(12n+7)^3}
\nonumber \\ &-&
\frac{9}{(12n+8)^3}+\frac{3}{(12n+10)^3}+\frac{3}{(12n+11)^3}+\frac{2}{(12n+12)^3} \Big)
\end{eqnarray}

Analogously we find for $\zeta(5)$:
\begin{eqnarray}
\zeta(5)&=&\frac{3\pi^2}{29}\zeta(3)-\frac{25\pi^4}{12528}+\frac{\pi^4}{1044}ln\left(\frac{\pi}{3}\right)
~-~\frac{144\pi^2}{29}\sum_{n=1}^{\infty}(-)^{n+1}c_n \left(\frac{\pi}{6\sqrt6}\right)^{2n}~, \nonumber \\
\end{eqnarray}
with
\begin{eqnarray}
c_1~=~\sum_{n=1}^{\infty}\frac{n(n+1)P^{(1)}(n)}{(2n-1)2n(2n+1)(2n+2)(2n+3)(2n+4)}\left(\frac{\pi}{6\sqrt6}\right)^{2n}~,
\end{eqnarray}
\begin{eqnarray}
c_2~=~\sum_{n=1}^{\infty}\frac{n(n+1)P^{(2)}(n+1)}{(2n+1)(2n+2)(2n+3)(2n+4)(2n+5)(2n+6)}\left(\frac{\pi}{6\sqrt6}\right)^{2n}~, \nonumber \\
\end{eqnarray}
\begin{eqnarray}
c_3~=~\sum_{n=1}^{\infty}\frac{n(n+1)P^{(3)}(n+2)}{(2n+3)(2n+4)(2n+5)(2n+6)(2n+7)(2n+8)}\left(\frac{\pi}{6\sqrt6}\right)^{2n}~, \nonumber \\
\end{eqnarray}
\begin{eqnarray}
c_4~=~\sum_{n=1}^{\infty}\frac{n(n+1)P^{(4)}(n+3)}{(2n+5)(2n+6)(2n+7)(2n+8)(2n+9)(2n+10)}\left(\frac{\pi}{6\sqrt6}\right)^{2n}~, \nonumber \\
\end{eqnarray}
Table IV shows the corresponding computational results:
\begin{table}
\begin{center}
\begin{tabular}[t]{||c||c|c|c||}
   \hline
   &&&\\
   $\zeta$(5) & $\zeta$(5)-Zeta series & $\zeta$(5)-(Zeta series+BP) & $\zeta$(5)-BBP formula \cite{kun13} \\ 
   &&&\\
   \hline
   1st order &&& \\
   n=1 & $\delta$=0.6*10$^{-06}$ &  $\delta$=0.1*10$^{-06}$ & $\delta$=0.7*10$^{-06}$ \\ 
   \hline
   2nd order &&& \\
   n=2 & $\delta$=0.4*10$^{-08}$ &  $\delta$=0.1*10$^{-09}$ & $\delta$=0.1*10$^{-08}$ \\ 
   \hline
   3rd order &&& \\
   n=3 & $\delta$=0.3*10$^{-10}$ &  $\delta$=0.1*10$^{-12}$ & $\delta$=0.8*10$^{-11}$ \\ 
   \hline
   4th order &&& \\
   n=4 & $\delta$=0.3*10$^{-12}$ &  $\delta$=0.1*10$^{-15}$ & $\delta$=0.1*10$^{-12}$ \\ 
   \hline
\end{tabular}
\caption{Approximate computation of $\zeta$(5) as a function of the summation index n by
use of Eq.~(1.15) without and with use of the polynomial representation. The numerical errors are compared
to the BBP-type formula \cite{kun13}}
\end{center}
\end{table}
where the numerical comparison has been performed between our approach and the BBP-formalism by use of the following
formula \cite{broad98,kun13}:
\begin{eqnarray}
\frac{31}{32}\zeta(5)&-&\frac{343}{99360}\pi^4ln(2)+\frac{5}{2484}\pi^2(ln(2))^3-\frac{2}{1035}(ln(2))^5
\nonumber \\ &=& \frac{128}{69}\sum_{n=1}^{\infty} \left(\frac{1}{\sqrt2}\right)^n \frac{cos(\frac{n\pi}{4})}{n^5}
-\frac{20}{69}\sum_{n=1}^{\infty} \left(\frac{1}{2}\right)^n \frac{1}{n^5}
\end{eqnarray}
 
For $\zeta(7)$ we found: 
\begin{eqnarray}
\zeta(7)&=& \frac{72\pi^2}{659}\zeta(5)-\frac{2\pi^4}{1977}\zeta(3)+\frac{49\pi^6}{5337900}
+\frac{\pi^6}{266895}ln\left(\frac{\pi}{3}\right)
\nonumber \\ &+&\frac{3456\pi^4}{5931}\sum_{n=1}^{\infty}(-)^{n+1}c_n
\left(\frac{\pi}{6\sqrt6}\right)^{2n}~, \nonumber \\ 
\end{eqnarray}
with
\begin{eqnarray}
c_1~&=&~\sum_{n=1}^{\infty}\frac{n(n+1)P^{(1)}(n)\left(\frac{\pi}{6\sqrt6}\right)^{2n}}{(2n-1)2n(2n+1)(2n+2)(2n+3)(2n+4)(2n+5)(2n+6)}
\end{eqnarray}
\begin{eqnarray}
c_2~&=&~\sum_{n=1}^{\infty}\frac{n(n+1)P^{(2)}(n+1)\left(\frac{\pi}{6\sqrt6}\right)^{2n}}{(2n+1)(2n+2)(2n+3)(2n+4)(2n+5)(2n+6)(2n+7)(2n+8)}
\nonumber \\
\end{eqnarray}
\begin{eqnarray}
c_3~&=&~\sum_{n=1}^{\infty}\frac{n(n+1)P^{(3)}(n+2)\left(\frac{\pi}{6\sqrt6}\right)^{2n}}{(2n+3)(2n+4)(2n+5)(2n+6)(2n+7)(2n+8)(2n+9)(2n+10)}
\nonumber \\
\end{eqnarray}
\begin{eqnarray}
c_4~&=&~\sum_{n=1}^{\infty}\frac{n(n+1)P^{(4)}(n+3)\left(\frac{\pi}{6\sqrt6}\right)^{2n}}{(2n+5)(2n+6)(2n+7)(2n+8)(2n+9)(2n+10)(2n+11)(2n+12)}
\nonumber \\
\end{eqnarray}
Table V shows the corresponding computational results:
\begin{table}
\begin{center}
\begin{tabular}[t]{||c||c|c|c||}
   \hline
   &&&\\
   $\zeta$(7) & $\zeta$(7)-Zeta series & $\zeta$(7)-(Zeta series+BP) & $\zeta$(7)-BBP formula $\;$ $\;$ $\;$  \\ 
   &&&\\
   \hline
   1st order &&& \\
   n=1 & $\delta$=0.7*10$^{-08}$ &  $\delta$=0.2*10$^{-08}$ &  \\ 
   \hline
   2nd order &&& \\
   n=2 & $\delta$=0.3*10$^{-10}$ &  $\delta$=0.1*10$^{-11}$ &  \\ 
   \hline
   3rd order &&& \\
   n=3 & $\delta$=0.2*10$^{-12}$ &  $\delta$=0.8*10$^{-15}$ &  \\ 
   \hline
   4th order &&& \\
   n=4 & $\delta$=0.2*10$^{-14}$ &  $\delta$=0.6*10$^{-18}$ &  \\ 
   \hline
\end{tabular}
\caption{Approximate computation of $\zeta$(7) as a function of the summation index n by
use of Eq.~(1.16) without and with use of the polynomial representation. For $\zeta$(7) no BBP-type formula is
available for reasons of comparison.}
\end{center}
\end{table}
where no BBP-type formula is available for $\zeta$(7).

\section{Approximations of zeta numbers using B$_{2n}$ and B$_{2n-2}$} 
We present a formula which allows to compute the Bernoulli number B$_{2n}$ as a function of the Bernoulli number B$_{2n-2}$ only.
This formula decouples B$_{2n}$ from all other Bernoulli numbers. We simply compute (3.32) for k=4 with the help of (3.3).
It follows:
\begin{eqnarray}
\frac{B_{2n}}{2n}~=~\frac{B_{2n-2}}{2n-2}~+~(-)^{n+1}\frac{(2n)!}{2^{2n}6^n}
               \left(\begin{array}{c} 
                               2n\\
                               4
			\end{array}\right)^{-1}\sum_{l=1}^{n}(-)^{l+1}
               \left(\begin{array}{c} 
                               n+4-l\\
                               4
			\end{array}\right) P^{(l)}(n)~. \nonumber \\
\end{eqnarray}
This formula allows the computation of $B_{2n}$ as a function of B$_{2n-2}$ only by use of the corresponding 
Bernoulli-type polynomials $P^{(l)}(n)$. Furthermore, it allows for an approximation of the odd zeta numbers
which works faster by about one order of magnitude when compared to (5.2).
It follows first:
\begin{eqnarray}
\zeta(2n)=\frac{(\zeta(2))^n}{2(2n-1)}\sum^{n+1}_{l=1}(-)^{l+1}
               \left(\begin{array}{c} 
                               n+5-l\\
                               4
			\end{array}\right) P^{(l)}(n+1)~
-\frac{2n(2n+1)}{4\pi^2}\zeta{(2n+2)}~. \nonumber \\
\end{eqnarray}
With this we can write:
\begin{eqnarray}
\zeta(3)=\frac{\pi^4}{15552}+\frac{\pi^2}{8}-\frac{\pi^2}{12}ln\left(\frac{\pi}{3}\right)+\frac{\pi^2}{6}
\sum^{\infty}_{n=1} (-)^{n+1} c_n \left(\frac{\pi}{6\sqrt6}\right)~,
\end{eqnarray}
with
\begin{eqnarray}
c_1~&=&~\frac{1}{24}\Big[\sum_{n=1}^{\infty}\frac{(n+1)(n+2)(n+3)(n+4)P^{(1)}(n+1)}{(2n-1)2n(2n+1)(2n+2)(2n+3)(2n+4)}
\nonumber \\ &=&
\frac{24n(n+1)P^{(1)}(n)}{(2n-1)2n(2n+1)(2n+2)(2n+3)(2n+4)} \Big] \left(\frac{\pi}{6\sqrt6}\right)^{2n}~,
\end{eqnarray}
\begin{eqnarray}
c_2~&=&~\frac{1}{24}\Big[\sum_{n=1}^{\infty}\frac{n(n+1)(n+2)(n+3)P^{(2)}(n+1)}
{(2n-1)2n(2n+1)(2n+2)(2n+3)(2n+4)} \nonumber \\ &=&
\frac{24n(n+1)P^{(2)}(n+1)\left(\frac{\pi}{6\sqrt6}\right)^{2}}{(2n+1)(2n+2)(2n+3)(2n+4)(2n+5)(2n+6)} \Big]
\left(\frac{\pi}{6\sqrt6}\right)^{2n}~,
\end{eqnarray}
\begin{eqnarray}
c_3~&=&~\frac{1}{24}\Big[\sum_{n=1}^{\infty}\frac{n(n+1)(n+2)(n+3)P^{(3)}(n+2)}
{(2n+1)(2n+2)(2n+3)(2n+4)(2n+5)(2n+6)} \nonumber \\ &=&
\frac{24n(n+1)P^{(3)}(n+2)\left(\frac{\pi}{6\sqrt6}\right)^{2}}{(2n+3)(2n+4)(2n+5)(2n+6)(2n+7)(2n+8)} \Big]
\left(\frac{\pi}{6\sqrt6}\right)^{2n}~,
\end{eqnarray}
\begin{eqnarray}
c_4~&=&~\frac{1}{24}\Big[\sum_{n=1}^{\infty}\frac{n(n+1)(n+2)(n+3)P^{(4)}(n+3)}
{(2n+3)(2n+4)(2n+5)(2n+6)(2n+7)(2n+8)} \nonumber \\ &=&
\frac{24n(n+1)P^{(4)}(n+3)\left(\frac{\pi}{6\sqrt6}\right)^{2}}{(2n+5)(2n+6)(2n+7)(2n+8)(2n+9)(2n+10)} \Big]
\left(\frac{\pi}{6\sqrt6}\right)^{2n}~,
\end{eqnarray}
and this way for higher coefficients c$_n$. Using these first four coefficients only we find an approximate value for
$\zeta(3)$ with an error $\delta <$ 0.4*10$^{-14}$. The convergence is about 4 orders of magnitude with the summation
index n. This procedure allows for a systematic improvement of the convergence behavior if higher values k=6 or 
k=8 are used in (3.2). As an example we present the corresponding formula for k=6:
\begin{eqnarray}
\left(\begin{array}{c} 
        2n \\
        6
\end{array}\right) 
\left(\frac{B_{2n}}{2n}-5\frac{B_{2n-2}}{2n-2}+4\frac{B_{2n-4}}{2n-4} \right)
(-)^{n+1}\frac{(2n)!}{2^{2n}6^n}\sum^n_{l=1}(-)^{l+1}
\left(\begin{array}{c} 
        n+6-l \\
        6
\end{array}\right)P^{(l)}(n)~. \nonumber \\ 
\end{eqnarray}
Our results establish a new and very fast option to compute zeta and related numbers. As an outlook, it could
be of great interest to combine our formalism with the BBP approach by use of corresponding polylogarithmic
identities to further improve the convergence properties in explicit computations.

\section{SUMMARY}
In summary, we have presented a unique computational scheme for the explicit calculation of the Riemann
$\zeta$ function and its first derivatives at all positive and negative integer values. This way we have
shown that all these numbers are directly attributed to Bernoulli numbers, but with an increasing level
of complexity when going, for example, from $\zeta$(2) to a related sub-sum like U$_4$. The computational
scheme is based on a new polynomial representation of the Bernoulli numbers in connection with Bendersky's
L-numbers, which appear in context with the logarithmic Gamma function. As a first application we performed
approximate calculations of $\zeta$(3), $\zeta$(5) and $\zeta$(7) in terms our polynomial representation,
where this computational procedure is applicable to all $\zeta$-values with integer arguments, as well as
to related numbers like Catalan's constant. Finally, we have shown that a computation of $B_{2n}$ as a
function of B$_{2n-2}$ or, for example, of B$_{2n-2}$ and B$_{2n-4}$ only is possible by use of the 
polynomial representation. The result is a further improved approximate computation of these numbers.

\section{APPENDIX A}
In analogy to example (3.5) for k=0 it results for k=1 and k=2:
\begin{eqnarray}
\sum^n_{l=1}(-)^{l+1}~l~P^{(n+1-l)}(n)~=~\frac{(2^{2n}-2)6^n}{2n(2n!)}B_{2n},
\end{eqnarray}
or
\begin{eqnarray}
\sum^n_{l=1}(-)^{l+1}~l~P^{(l)}(n)~=~(-)^{n+1}\frac{(4n-(n-1)2^{2n})6^n}{2n(2n!)}B_{2n}~,
\end{eqnarray}
with
\begin{eqnarray}
\sum^n_{l=1}(-)^{l+1}~l~P^{(n+1-l)}(n)~= \hspace*{4.2cm}\\  \nonumber
~(-)^{n+1}\frac{2n(2^{2n}-2)}{4n-(n-1)2^{2n}}
\sum^n_{l=1}(-)^{l+1}~l~P^{(l)}(n)~.
\end{eqnarray}
For k=2 we find
\begin{eqnarray}
\sum^n_{l=1}(-)^{l+1}~l^2~P^{(n+1-l)}(n)~=~\frac{(2+(2n-2)2^{2n})6^n}{(2n!)}B_{2n},
\end{eqnarray}
or
\begin{eqnarray}
\sum^n_{l=1}(-)^{l+1}~l^2~P^{(l)}(n)~=\hspace*{4.2cm}\\  \nonumber
~(-)^{n+1}\frac{(4n(2n+3)+(n^2-6n+1)2^{2n})6^n}{(2n!)2n}B_{2n},
\end{eqnarray}
with
\begin{eqnarray}
\sum^n_{l=1}(-)^{l+1} l^2~P^{(n-l+1)}(n)~= \hspace*{4.2cm}\\  \nonumber
~(-)^{n+1}~\frac{2n(2+(2n-2)2^{2n})}{4n(2n+3)+(n^2-6n+1)2^{2n}} \sum^n_{l=1}(-)^{l+1} l^2~P^{(l)}(n)~,
\end{eqnarray}
where the simple symmetry property 
\begin{eqnarray}
\sum^n_{l=1}(-)^{l+1}l^m~ P^{(l)}(n)&=& \nonumber \\ &&
(-)^{n+1}\sum^n_{l=1}(-)^{l+1} (n+1-l)^m~P^{(n+1-l)}(n)~,
\end{eqnarray}
which we have used before in the case of m=1 has been applied here again for m=0,1 and 2.

\section{APPENDIX B}
Here we show the Bernoulli polynomials up to n=6. It follows:
\begin{eqnarray}
P^{(1)}(n)~=~\frac{1}{n}
\end{eqnarray}
\begin{eqnarray}
P^{(2)}(n)~=~\frac{3}{2*5}
\end{eqnarray}
\begin{eqnarray}
P^{(3)}(n)~=~\frac{3(21n-43)}{2^3*5^2*7}
\end{eqnarray}
\begin{eqnarray}
P^{(4)}(n)=\frac{63n^2-387n+590}{2^4*5^3*7}
\end{eqnarray}
\begin{eqnarray}
P^{(5)}(n)=\frac{3(4851n^3-59598n^2+242737n+327210)}{2^75*^4*7^2*11}
\end{eqnarray}
\begin{eqnarray}
P^{(6)}(n)=\frac{3(189189n^4-3873870n^3+29616015n^2-100104550n+126087736)}{2^8*5^6*7^2*11*13} \nonumber \\
\end{eqnarray}

\makeatletter
\renewcommand{\@biblabel}[1]{[#1]\hfill}
\makeatother

\end{document}